\documentclass[12pt]{amsart}
\usepackage{amssymb}
\usepackage{amsmath}
\newtheorem{Lemma}{Lemma}[section]

\newtheorem{Theorem}[Lemma]{Theorem}
\newtheorem{Corollary}[Lemma]{Corollary}
\newtheorem{Definition}[Lemma]{Definition}
\newtheorem{Proposition}[Lemma]{Proposition}

\def\C{\mathbb{C}}
\def\Cas{\mathrm{Cas}}
\def\Der{\mathrm{Der}}
\def\e{\varepsilon}
\def\ext{\mathrm{ext}}
\def\End{\mathrm{End}}
\def\g{\mathfrak{g}}
\def\G{\mathbf{G}}

\def\H{\mathbb{H}}
\def\Hodge{\mathrm{Hodge}}
\def\hol{\mathfrak{hol}}
\def\Hol{\mathrm{Hol}}
\def\Hom{\mathrm{Hom}}
\def\id{\mathrm{id}}
\def\ins{\lrcorner}
\def\k{\mathfrak{k}}
\def\ker{\mathrm{ker}}
\def\L{\Lambda}
\def\op{\diamond}
\def\pr{\mathrm{pr}}
\def\R{\mathbb{R}}
\def\Ric{\mathrm{Ric}}
\def\S{\mathrm{Sym}}
\def\so{\mathfrak{so}}
\def\SO{\mathrm{SO}}

\def\Sp{\mathrm{Sp}}

\def\Spin{\mathrm{Spin}}
\def\su{\mathfrak{su}}
\def\SU{\mathbf{SU}}

\def\X{\mathfrak{X}}
\title{The Standard Laplace Operator}
\author{Uwe Semmelmann \& Gregor Weingart}
\address{Uwe Semmelmann\\
Institut f\"ur Geometrie und Topologie \\
Fachbereich Mathematik\\
Universit\"at Stuttgart\\
Pfaffenwaldring 57 \\
70569 Stuttgart, Germany}
\email{uwe.semmelmann@mathematik.uni-stuttgart.de}
\address{Gregor Weingart\\
Instituto de Matem\'aticas (Unidad Cuernavaca)\\
Universidad Nacional Aut\'onoma de Mexico\\
Avenida Universidad s/n\\
Colonia Lomas de Chamilpa\\
62210 Cuernavaca, Morelos, Mexico}
\email{gw@matcuer.unam.mx}
\date{\today}
\begin{document} 
\begin{abstract}
 The standard Laplace operator is a generalization of the Hodge Laplace
 operator on differential forms to arbitrary geometric vector bundles,
 alternatively it can be seen as generalization of the Casimir operator
 acting on sections of homogeneous vector bundles over symmetric spaces
 to general Riemannian manifolds. Stressing the functorial aspects of the
 standard Laplace operator $\Delta$ with respect to the category of geometric
 vector bundles we show that the standard Laplace operator commutes not only
 with all homomorphisms, but also with a large class of natural first order
 differential operators between geometric vector bundles. Several examples
 are included to highlight the conclusions of this article.

 \medskip
 \begin{center}
  \parbox{300pt}{\textit{MSC 2010:}\ 53C21; 53C26, 58A14.}
  \\[5pt]
  \parbox{300pt}{\textit{Keywords:}\ Laplace operator,
  Riemannian manifold, special holonomy.}
 \end{center}
\end{abstract}
\maketitle
\section{Introduction}
 \noindent
 The standard Laplace operator on geometric vector bundles can be seen
 as the quintessence of and motivation for the work of the authors on
 special holonomy in general and Weitzenb\"ock formulas in particular.
 Its mere existence is some kind of miracle, because it generalizes both
 the Casimir operator on symmetric spaces to general Riemannian manifolds
 and the Hodge Laplace operator on differential forms to more general
 vector bundles.

 Underlying the construction of the standard Laplace operator is the concept
 of geometric vector bundles on a Riemannian manifold $(\,M,\,g\,)$
 well adapted to the geometry arising from a reduction of its holonomy group.
 Intuitively speaking a geometric vector bundle is a vector bundle on $M$
 associated to the principal subbundle $\Hol\,M\,\subset\,\SO\,M$ of the
 bundle of oriented orthonormal frames encoding this holonomy reduction,
 a precise axiomatic definition avoiding all the subtleties of holonomy
 reduction and principal bundles is given below. By construction geometric
 vector bundles form a category of vector bundles with connection on the
 manifold $M$. The characteristic properties of geometric vector bundles
 among all vector bundles with connection allows us however to construct
 a curvature endomorphism
 $$
  q(\,R\,)\;\;\in\;\;\Gamma(\;\End\,VM\;)
 $$
 for every geometric vector bundle $VM$ parametrized by the curvature tensor
 $R$ of a metric connection $\bar\nabla$ on the tangent bundle of $M$, which
 commutes with all homomorphisms of geometric vector bundles. Needless to say
 this curvature endomorphism has been considered in numerous special cases
 before and we will recall several interesting examples below, however we
 believe that only the category of geometric vector bundles brings the
 construction of the curvature endomorphism $q(R)$ to the point. The
 standard Laplace operator associated to a geometric vector bundle is
 the sum of the rough Laplacian with the curvature endomorphism:
 \begin{equation}\label{standard}
  \Delta\;\;:=\;\;\nabla^*\nabla\;+\;q(\,R\,)\ .
 \end{equation}
 By construction $\Delta$ is a Laplace type operator which commutes with all
 homomorphisms of geometric vector bundles. Of course the same statement
 is true for every linear combination of $\nabla^*\nabla$ and $q(R)$, the
 variation of the curvature tensor $R$ under the Ricci flow $\delta g\,=\,
 2\Ric$ for example can be written neatly in the form $\delta R\,=\,
 (\,\nabla^*\nabla\,+\,\frac12\,q(R)\,)\,R\,+\,\frac12\,\Der_\Ric R$.
 The decisive advantage of the linear combination (\ref{standard}) is
 however that $\Delta$ commutes not only with all homomorphisms of geometric
 vector bundles, but also with many natural first order differential
 operators $D$. More precisely the commutator is calculated in Theorem
 \ref{theorem}
 \begin{equation}\label{comm}
  [\;\Delta,\;D\;]
  \;=\;\;
  -\,\sigma_\op(\;q(\,\bar\nabla R\,)\;-\;\delta R\;) ,
 \end{equation}
 and depends only on the principal symbol $\sigma_\op$ of the natural first
 order differential operator $D$ and the covariant derivative $\bar\nabla R$
 of the curvature tensor $R$ of the given metric connection $\bar\nabla$ on
 $TM$. On a symmetric space for example $\bar\nabla R\,=\,0$ and thus $\Delta$
 commutes with all natural differential operators between geometric vector
 bundles, showing a priori that it agrees with the Casimir operator on
 homogeneous vector bundles.

 Illustrating the intrinsic beauty of the construction of the standard
 Laplace operator in many explicit and rather well--known examples makes
 this article somewhat reminiscent of a review article, nevertheless it
 contains a couple of new ideas besides establishing the general commutator
 formula (\ref{comm}) and discussing its applications; both the underlying
 concept of geometric vector bundles and the presentation of the curvature
 endomorphism $q(R)$ as an integral over sectional curvatures certainly
 deserve a more detailed analysis in the future.

 \medskip\noindent
 In the second section of this article we introduce the concept of geometric
 vector bundles and discuss several constructions of the curvature
 endomorphism $q(R)$, in particular its integral presentation. In the
 third section we prove the commutator formula (\ref{comm}) in Theorem
 \ref{theorem} and Corollary \ref{corollary} and in the fourth section
 we sketch some direct applications of this commutator formula to the
 geometry of normal homogeneous spaces, nearly K\"ahler manifolds in
 dimension six, manifolds with $\G_2$-- and $\Spin(7)$--holonomy and
 quaternion--K\"ahler manifolds.
\section{Preliminaries}\label{prelim}
\subsection{Geometric Vector Bundles}
 In differential geometry it is well--known that many interesting flavors
 of geometry come along with a corresponding reduction of the holonomy
 group, in fact this observation is one of the basic tenets of Cartan's
 generalization of the Erlangen program nowadays called Cartan geometry.
 In Riemannian geometry such a holonomy reduction can be thought of as a
 subbundle $\Hol\,M\,\subset\,\SO\,M$ of the orthonormal frame bundle of
 a manifold $M$ which is a principal bundle itself under the induced action
 of a closed subgroup $\Hol\,\subset\,\SO\,T$ and tangent to a metric
 connection $\bar\nabla$ on the tangent bundle $TM$. The most direct
 implication of such a holonomy reduction is that the metric connection
 $\bar\nabla$ induces a principal connection on $\Hol\,M$. The association
 functor, which turns a representation $V$ of $\Hol$ into the associated
 vector bundle $\Hol\,M\times_\Hol V$ on $M$, thus becomes a functor
 $\mathrm{Asc}_{\Hol\,M}:\,\mathrm{Rep}_\Hol\,\rightsquigarrow
 \,\mathrm{Vect}^\nabla_M$ from the category of representations of $\Hol$
 to the category of vector bundles with connection over $M$ under parallel
 homomorphisms.

 In essence a geometric vector bundle on $M$ is a vector bundle in the
 image of this association functor. In order to avoid the subtleties of
 holonomy reductions and the association functor we will provide a more
 direct axiomatic definition of a geometric vector bundle though, which
 characterizes a geometric vector bundle $VM$ by an infinitesimal
 representation of the holonomy algebra subbundle $\hol\,M\,\subset\,
 \so\,TM\,\cong\,\L^2TM$ for the given metric connection $\bar\nabla$:

 \begin{Definition}[Geometric Vector Bundles]
 \hfill\label{gvb}\break
  A geometric vector bundle on a Riemannian manifold $(\,M,\,g\,)$ endowed
  with a metric connection $\bar\nabla$ on its tangent bundle is a vector
  bundle $VM$ endowed with a connection $\nabla$ and a parallel infinitesimal
  representation $\star:\,\hol\,M\otimes VM\longrightarrow VM$ of the holonomy
  algebra bundle $\hol\,M\,\subset\,\L^2TM$ such that the curvature $R^\nabla$
  of the connection $\nabla$ is determined via
  $$
   R^\nabla_{X,\,Y}\;\;=\;\;R_{X,\,Y}\,\star
  $$
  for all $X,\,Y\,\in\,\Gamma(\,TM\,)$ by the curvature tensor $R$ of the
  metric connection $\bar\nabla$ on $TM$.
 \end{Definition}

 \noindent
 Needless to say geometric vector bundles form a category, the morphisms
 in this category are parallel homomorphisms $F:\,VM\longrightarrow\tilde VM$
 of vector bundles with connections which commute with the infinitesimal
 representation of the holonomy algebra bundle $\hol\,M$. Evidently this
 category is closed under taking duals, direct sums, tensor products and
 exterior and symmetric powers, to name only a few constructions of linear
 algebra. The relevance of geometric vector bundles for this article is that
 the curvature endomorphism and the standard Laplace operator constructed
 below are essentially functors in the sense that they commute with all
 morphisms in the category of geometric vector bundles.

 \medskip
 {\bf Example 1:}
 A classical example of a geometry related to a holonomy reduction is
 K\"ahler geometry: On a K\"ahler manifold $(\,M,\,g,\,J\,)$ the bundle
 of orthonormal frames reduces to complex linear orthonormal frames with
 its unitary structure group $\mathbf{U}(\,n\,)\,\cong\,\Hol\subset\,\SO\,T$.
 In turn the holonomy algebra bundle equals the bundle $\hol\,M\,\subset\,
 \so\,TM$ of skew symmetric endomorphisms of $TM$ commuting with $J$ with
 projection $\pr_\hol:\,\L^2TM\longrightarrow\hol\,M$ given by:
 $$
  \pr_\hol(\,X\,\wedge\,Y\,)
  \;\;=\;\;
  \pr^{1,0}X\,\wedge\,\pr^{0,1}Y\;+\;\pr^{0,1}X\,\wedge\,\pr^{1,0}Y\ .
 $$
 It is well known in this case that the geometric vector bundle of
 antiholomorphic differential forms $VM\,=\,\L^{0,\bullet}T^*M$ is a
 Clifford bundle with parallel Clifford multiplication defined as
 $$
  X\,\bullet
  \;\;:=\;\;
  \sqrt{2}\,(\;\pr^{1,0}X^\flat\,\wedge\;-\;\pr^{0,1}X\,\ins\;)\ ,
 $$
 moreover its fiber $V_xM$ is a spinor module for the Clifford algebra
 $\mathrm{Cl}(\,T_xM,\,g_x\,)$ in every point $x\,\in\,M$. Nevertheless
 the vector bundle $VM$ is not considered to be the spinor bundle $\Sigma M$
 of a K\"ahler manifold $M$, and K\"ahler manifolds as simple as $\C P^2$
 are not even spin.  The concept of geometric vector bundles clarifies this
 apparent contradiction: The infinitesimal representation
 $\star:\,\hol\,M\otimes VM\longrightarrow VM$ of the geometric vector
 bundle $VM$ is the restriction of the standard representation of $\so\,TM$
 on differential forms to $\hol\,M$ and $VM$, for the spinor bundle $\Sigma M$
 however the infinitesimal representation $\star_\Sigma:\,\so\,TM\otimes
 \Sigma M\longrightarrow\Sigma M$ is declared by way of axiom to be induced
 by Clifford multiplication with bivectors:
 $$
  (\,X\,\wedge\,Y\,)\,\star_\Sigma
  \;\;:=\;\;
  \tfrac12\,(\;X\,\bullet\,Y\,\bullet\;+\;g(\,X,\,Y\,)\,\id_{\Sigma M}\;)\ .
 $$
 Of course we may simply redefine the infinitesimal representation of
 $VM$ to be the restriction of $\star_\Sigma$ to the actual holonomy bundle
 $\hol\,M\,\subset\,\so\,TM$ of K\"ahler geometry. This replacement however
 will not result in a geometric vector bundle $(\,VM,\,\nabla,\,\star_\Sigma
 \,)$ isomorphic to the spinor bundle $\Sigma M$ as a geometric vector bundle,
 because the difference between the representations
 $$
  \pr_\hol(\,X\,\wedge\,Y\,)\,\star_\Sigma
  \;-\;\pr_\hol(\,X\,\wedge\,Y\,)\,\star
  \;\;=\;\;
  -\,i\,g(\,JX,\,Y\,)\,\id_{\Sigma M}
 $$
 implies that the action of the curvature tensor $R_{X,\,Y}\star_\Sigma
 \,=\,R_{X,\,Y}\star\,+\,i\,\Ric(JX,Y)\,\id_{\Sigma M}$ under $\star_\Sigma$
 does not agree with the curvature $R_{X,\,Y}\,=\,R_{X,\,Y}\star$ of the
 given connection $\nabla$ on $VM$ unless the K\"ahler manifold $M$ is
 actually Ricci flat and thus a Calabi--Yau manifold.
\subsection{Generalized Gradients}
 On sections of a geometric vector bundle $VM$ we have a distinguished
 set of natural first order differential operators called generalized
 gradients defined by projecting the covariant derivative $\nabla$ to
 its isotypical components. More precisely we consider the canonical
 decomposition $T\otimes V\,=\,\bigoplus V_\e$ of the $\Hol$--representation
 $T\otimes V$ into isotypical components with corresponding projections
 $\pr_\e:\,T\otimes V\longrightarrow V_\e\,\subset\,T\otimes V$, which
 become parallel projections $\pr_\e:\,TM\otimes VM\longrightarrow V_\e M$
 between the corresponding geometric vector bundles. The {\em generalized
 gradient} $P_\e$ is defined for all $\e$ as the composition:
 $$
  P_\e:\;\;\Gamma(\,VM\,)
  \;\stackrel\nabla\longrightarrow\;
  \Gamma(\,T^*M\,\otimes\,VM\,)
  \;\stackrel\cong\longrightarrow\;
  \Gamma(\,TM\,\otimes\,VM\,)
  \;\stackrel{\pr_\e}\longrightarrow\;
  \Gamma(\,V_\e M\,)\ .
 $$
 In a similar vein we may define natural second order differential
 operators on sections of $VM$ by taking constant linear combinations
 of the operators $P^*_\e P^{\phantom{*}}_\e$. Certainly the most
 important example of this kind of second order differential operators
 is the so called {\em rough Laplacian}
 $$
  \nabla^* \nabla 
  \;\;:=\;\;
  -\,\sum_\mu\Big(\,\nabla_{E_\mu}\,\nabla_{E_\mu}\;-\;
  \nabla_{\bar\nabla_{E_\mu}E_\mu}\,\Big)
  \;\;=\;\;
  \sum_\e P^*_\e\,P^{\phantom{*}}_\e\ ,
 $$
 where the sum is over some local orthonormal frame $\{\,E_\mu\,\}$ of the
 tangent bundle $TM$. Besides generalized gradients we will consider
 natural first order differential operators arising from an arbitrary
 $\Hol$--equivariant homomorphism $\sigma_\op:\,T\otimes V\longrightarrow
 \tilde V$ for two representations $V$ and $\tilde V$ as well, if only to
 simplify notation. More precisely we compose $\nabla$ with the parallel
 extension of $\sigma_\op$ to vector bundles to define the natural first
 order differential operator:
 $$
  D_\op:\;\;\Gamma(\,VM\,)
  \;\stackrel\nabla\longrightarrow\;
  \Gamma(\,T^*M\,\otimes\,VM\,)
  \;\stackrel\cong\longrightarrow\;
  \Gamma(\,TM\,\otimes\,VM\,)
  \;\stackrel{\sigma_\op}\longrightarrow\;
  \Gamma(\,\tilde VM\,)\ .
 $$
 Alternatively we may write $D_\op$ as a sum over a local orthonormal basis
 $\{\,E_\mu\,\}$ of the tangent bundle with respect to the Riemannian metric
 $g$ in order to make this definition resemble the definition of the Dirac
 operator on Clifford bundles with the Clifford multiplication $\bullet$
 replaced by the parallel bilinear operation $\sigma_\op:\,TM\otimes VM
 \longrightarrow\tilde VM,\,X\otimes v\longmapsto X\,\op\,v$:
 \begin{equation}\label{dsum}
  D_\op\;\;:=\;\;\sum_\mu E_\mu\,\op\,\nabla_{E_\mu}\ .
 \end{equation}
\subsection{Metric Connections with Torsion}
\label{mct}
 Working with metric, but not necessarily torsion free connections
 on Riemannian manifolds requires some care even for specialists in
 Riemannian geometry to ensure that torsion freeness is not assumed
 implicitly in some innocuously looking argument or other. Throughout
 this article we will consider metric connections $\bar\nabla$ with
 parallel skew--symmetric torsion only, their curvature tensors
 share most of the classical symmetries of the curvature tensor
 of the Levi--Civita connection. Working through a proof of the
 first Bianchi identity for torsion free connections we easily see
 that it extends to arbitrary connections $\bar\nabla$ on $TM$ with
 torsion tensor $T(X,Y)\,:=\,\bar\nabla_XY\,-\,\bar\nabla_YX\,-\,[X,Y]$:
 \begin{eqnarray*}
  R_{X,\,Y}Z\,+\,R_{Y,\,Z}X\,+\,R_{Z,\,X}Y
  &=&
  (\bar\nabla_XT)(Y,Z)\,+\,(\bar\nabla_YT)(Z,X)
  \,+\,(\bar\nabla_ZT)(X,Y)
  \\
  &&
  \quad+\,T(T(X,Y),Z)\,+\,T(T(Y,Z),X)\,+\,T(T(Z,X),Y)\ .
 \end{eqnarray*}
 In order to derive some meaningful conclusions from this rather tautological
 version of the first Bianchi identity let us first restrict the class of
 connections considered in this article:

 \begin{Definition}[Connections with Parallel Skew--Symmetric Torsion]
 \hfill\label{cpt}\break
  A metric connection $\bar\nabla$ on the tangent bundle $TM$ of a
  Riemannian manifold $(\,M,\,g\,)$ with torsion tensor $T$ is called a
  metric connection with parallel skew--symmetric torsion provided
  the expression $\theta(X,Y,Z)\,:=\,g(T(X,Y),Z)$ defines a parallel
  $3$--form $\theta\,\in\,\Gamma(\,\L^3T^*M\,)$.
 \end{Definition}

 \noindent
 Evidently the torsion tensor itself is parallel $\bar\nabla T\,=\,0$ for
 a metric connection $\bar\nabla$ with parallel skew--symmetric
 torsion, moreover the endomorphism $Y\longmapsto T(X,Y)$ is skew--symmetric
 with respect to the Riemannian metric $g$ for every $X\,\in\,\Gamma(TM)$.
 Due to $\bar\nabla T\,=\,0$ the tautological version of the first Bianchi
 identity reduces to the identity
 \begin{equation}\label{1bi}
  R(X,Y,Z,W)\;+\;R(X,Y,Z,W)\;+\;R(X,Y,Z,W)
  \;\;=\;\;
  \tfrac12\,g(T\wedge T)(X,Y,Z,W)
 \end{equation}
 valid for the curvature tensor $R(X,Y,Z,W)\,:=\,g(R_{X,Y}Z,W)$ of a
 connection $\bar\nabla$ with parallel skew--symmetric torsion,
 where $\frac12\,g(T\wedge T)\,\in\,\Gamma(\L^4T^*M)$ denotes the parallel
 $4$--form:
 \begin{eqnarray*}
  \lefteqn{\tfrac12\,g(\,T\,\wedge\,T\,)(\,X,\,Y,\,Z,\,W\,)}
  \qquad
  &&
  \\
  &:=&
  g(\,T(X,Y),\,T(Z,W)\,)\;+\;g(\,T(Y,Z),\,T(X,W)\,)\;+\;g(\,T(Z,X),\,T(Y,W)\,)
 \end{eqnarray*}
 In light of the modified first Bianchi identity (\ref{1bi}) a standard
 proof for the well--known symmetry $R(X,Y,Z,W)\,=\,R(Z,W,X,Y)$ of the
 curvature tensor of the Levi--Civita connection goes through more or less
 verbatim to prove $R(X,Y,Z,W)\,=\,R(Z,W,X,Y)$ for the curvature tensor
 $R$ of a metric connection $\bar\nabla$ with parallel skew--symmetric
 torsion (cf.~\cite{ivanov}). In particular the Ricci tensor $\Ric(X,Y)
 \,:=\,\sum_\mu R(X,E_\mu,E_\mu,Y)$ is symmetric as well.
\subsection{The Curvature Endomorphism}
 The standard identification of the Lie algebra bundle $\so\,TM$ of skew
 symmetric endomorphism on a Riemannian manifold with the bivector bundle
 $\so\,TM\,=\,\L^2TM$ realizes the holonomy bundle $\hol\,M$ as a subbundle
 of $\L^2TM$ endowed with a scalar product induced by the Riemannian metric
 on $M$. The parallel orthogonal projection map $\pr_\hol:\,\L^2TM
 \longrightarrow\hol\,M\,\subset\,\L^2TM$ to the holonomy subbundle
 allows us to define the standard curvature endomorphism for every
 geometric vector bundle:

 \begin{Definition}[Curvature Endomorphism]
 \hfill\label{ce}\break
  Let $\{\,E_\mu\,\}$ be a local orthonormal frame of the tangent bundle
  $TM$. The curvature endomorphism $q(R)\,\in\,\End\,VM$ is defined
  for every geometric vector bundle $VM$ as the sum:
  $$
   q(\,R\,)
   \;\;:=\;\;
   \tfrac12\,\sum_{\mu\nu}
   \pr_\hol(\,E_\mu\,\wedge\,E_\nu\,)\,\star\,R_{E_\mu,\,E_\nu}\ .
  $$
 \end{Definition}

 \noindent
 Of course the definition of $q(R)$ is independent of the choice $\{\,E_\mu
 \,\}$ of a local orthonormal frame for the tangent bundle $TM$. It is more
 important to observe however that this standard argument of linear algebra
 can be applied as well to the induced local orthonormal frame $\{\,E_\mu
 \wedge E_\nu\,\}$ of the Lie algebra bundle $\so\,TM\,=\,\L^2TM$. Considering
 the curvature tensor as an operator $R:\,\L^2TM\longrightarrow\hol\,M,\,
 X\wedge Y\longmapsto R_{X,\,Y},$ of vector bundles we may write
 $$
  q(\,R\,)
  \;\;=\;\;
  \tfrac12\,\sum_{\mu\nu}
  \pr_\hol(\,E_\mu\,\wedge\,E_\nu\,)\,\star\,R_{E_\mu,\,E_\nu}
  \;\;\stackrel!=\;\;
  \sum_\alpha(\,\pr_\hol\X_\alpha\,)\,\star\,R(\,\X_\alpha\,)\,\star
 $$
 with an arbitrary local orthonormal basis $\{\,\X_\alpha\,\}$ of $\so\,TM$.
 For a local orthonormal basis $\{\,\X_\alpha,\,\X^\perp_\beta\,\}$ adapted
 to the decomposition $\so\,TM\,=\,\hol\,M\oplus\hol^\perp M$ we find
 in particular
 $$
  q(\,R\,)
  \;\;=\;\;
  \sum_\alpha\X_\alpha\,\star\,R(\,\X_\alpha\,)\,\star\ ,
 $$
 because the sum over $\{\,\X^\perp_\beta\,\}$ vanishes. In particular
 $q(R)$ is a symmetric endomorphism for every geometric vector bundle and
 for every metric connection $\bar\nabla$ on the tangent bundle $TM$ with
 parallel skew--symmetric torsion, because $R\,\in\,\S^2\hol\,M\,\subset\,
 \S^2\L^2TM$ is symmetric for such connections according to Subsection
 \ref{mct}. In particular we can choose the local orthonormal basis
 $\{\,\X_\alpha\,\}$ of $\hol\,M$ to be a basis of eigenvectors of the
 curvature operator $R:\,\L^2TM\longrightarrow\hol\,M$, hence $q(R)\,\geq\,0$
 is a non--negative operator provided all eigenvalues of the curvature operator
 are non--positive and vice versa.

 A similar argument shows that $q(R)$ is hereditary under successive holonomy
 reductions in the sense that every geometric vector bundle $VM$ adapted to
 a holonomy reduction $\hol\,M$ remains a geometric vector bundle under a
 further reduction $\overline\hol\,M\,\subset\,\hol\,M$ of the holonomy algebra
 bundle; in this situation the curvature endomorphism $q(R)$ does not depend
 on which projection $\pr_\hol$ or $\pr_{\overline\hol}$ is chosen in
 Definition \ref{ce}. The projection $\pr_\hol$ for example makes no difference
 at all for a geometric vector bundle $VM$, whose infinitesimal representation
 $\star:\,\hol\,M\otimes VM\longrightarrow VM$ is actually the restriction of
 an infinitesimal representation of the generic holonomy bundle $\so\,TM\,
 \supset\,\hol\,M$. In particular the curvature endomorphism $q(R)$ equals
 the Ricci endomorphism for every geometric subbundle of the tangent bundle:
 \begin{eqnarray*}
  q(\,R\,)\,X
  &=&
  \tfrac12\,\sum_{\mu\nu}\pr_\hol(\,E_\mu\wedge E_\nu\,)\,\star\,
  R_{E_\mu,\,E_\nu}X
  \;\;\stackrel!=\;\;
  \tfrac12\,\sum_{\mu\nu}(\,E_\mu\wedge E_\nu\,)\,\star\,R_{E_\mu,\,E_\nu}X
  \\
  &=&
  \tfrac12\,\sum_{\mu\nu}\Big(\,g(E_\mu,R_{E_\mu,E_\nu}X)\,E_\nu
  \,-\,g(E_\nu,R_{E_\mu,E_\nu}X)\,E_\mu\,\Big)
  \;\;=\;\;
  \sum_\nu\Ric(\,E_\nu,\,X\,)\,E_\nu\ .
 \end{eqnarray*}

 \medskip
 {\bf Example 1:} The example motivating Definition \ref{ce} is certainly
 the forms representation $V\,=\,\L^\bullet T^*$ in generic holonomy
 $\hol\,=\,\so\,T$, in which a bivector $X \wedge Y\,\in\,\L^2T\,\cong\,
 \so\,T$ acts by $(X \wedge Y)\star\,:=\,Y^\flat\wedge X\ins\,-\,X^\flat
 \wedge Y\ins$. In this example the curvature endomorphism becomes:
 $$
  q(\,R\,)
  \;\;=\;\;
  \tfrac12\,\sum_{\mu\nu}
  (\,E_\mu\,\wedge\,E_\nu)\,\star\,R_{E_\mu,\,E_\nu}
  \;\;=\;\;
  -\,\sum_{\mu\nu}E_\mu^\flat\wedge\,E_\nu\,\ins\,R_{E_\mu,\,E_\nu}\ .
 $$
 The curvature term on the right hand side is well--known, because it
 appears in the classical Weitzenb\"ock formula for the Hodge--Laplace
 operator $\Delta_\Hodge\,:=\,dd^*\,+\,d^*d$ acting on $p$--forms, the
 classical Weitzenb\"ock formula $\Delta_\Hodge\,=\,\nabla^*\nabla\,+\,
 q(R)$ thus forms the blue print for the definition (\ref{standard}) of
 the standard Laplace operator $\Delta$. As a direct consequence of this
 definition the standard Laplace operator $\Delta\,=\,\Delta_\Hodge$
 equals the Hodge--Laplace operator for every geometric vector bundle
 occurring as a subbundle in the bundle of differential forms.

 \medskip
 {\bf Example 2:} In the generic case $\hol\,=\,\so\;T$ we consider the
 representation $V\,=\,\S^2_0T^*$ of trace--free symmetric $2$--tensors.
 In the notation of \cite{besse} the curvature endomorphism $q(R)$ acting
 on sections of the corresponding geometric vector bundle $\S^2_0T^*M$
 becomes the sum
 $$
  q(\,R\,)
  \;\;=\;\;
  2\,\raise10pt\hbox to0pt{$\;\scriptscriptstyle\circ$\hss}R
  \;+\;\Der_\Ric\ ,
 $$
 where $(\Der_\Ric h)(X,Y)\,:=\,h(\Ric\,X,Y)\,+\,h(X,\Ric\,Y)$ denotes
 the standard derivative extension of the Ricci endomorphism to bilinear
 forms, whereas the curvature operator $\raise10pt\hbox to0pt{$\;
 \scriptscriptstyle\circ$\hss}R$ is defined specifically for symmetric
 $2$--tensors $h\,\in\,\Gamma(\S^2T^*M)$ by:
 $$
  (\raise10pt\hbox to0pt{$\;\scriptscriptstyle\circ$\hss}R\,h)(\;X,\,Y\;)
  \;\;:=\;\;
  \sum_\mu h(\;R_{X,\,E_\mu}Y,\,E_\mu\;)\ .
 $$
 Given this ad hoc definition it is not even clear that $\raise10pt
 \hbox to0pt{$\;\scriptscriptstyle\circ$\hss}Rh$ is actually symmetric,
 although this can be shown by using the first Bianchi identity and the
 skew--symmetry of the endomorphism $Z\longmapsto R_{X,Y}Z$ for every
 metric connection on the tangent bundle $TM$.

 \medskip
 {\bf Example 3:}
 Consider a Riemann symmetric space $M\,=\,G/K$. The isometry group $G$
 of $M$ can be thought of as a principal bundle over $M$ with structure
 group structure equal to the isotropy group $K$ of the base point $eK\,
 \in\,M$ and it is well--known that this principal bundle agrees with the
 holonomy reduction $G\,\cong\,\Hol\,M$, in consequence geometric vector
 bundles on $M$ are homogeneous vector bundles and vice versa. There exists
 a unique non--degenerate invariant scalar product $g^\ext$ on the Lie
 algebra $\g$ of $G$, which restricts to the Riemannian metric $g$ on the
 symmetric complement $\mathfrak{p}\,\cong\,T_{eK}M\,\subset\,\g$ of the
 isotropy subalgebra $\k\,\subset\,\g$. The curvature endomorphism $q(R)$
 agrees on every geometric vector bundle with the Casimir operator of
 $\k$ defined as a sum over an orthonormal basis $\{\,\X_\alpha\,\}$
 of $\k$ with respect to $g^\ext$
 $$
  q(R)
  \;\;=\;\;
  \Cas^\k_{g^\ext}
  \;\;:=\;\;
  -\,\sum_\alpha\X_\alpha\star\,\X_\alpha\star
 $$
 (cf. Lemma 5.2 in \cite{au10}).

 \medskip
 {\bf Example 4:} Let $(\,M,\,g\,)$ be a spin manifold with the spin
 structure defined by a principal $\Spin\,T$--bundle $\Spin\,M$ lifting
 the bundle $\SO\,M$ of oriented orthonormal frames. For the geometric vector
 bundle $\Sigma M$ corresponding to the spinor representation $\Sigma$ of
 $\Spin\,T$ the curvature endomorphism acts as multiplication with the
 scalar curvature $\kappa\,\in\,C^\infty(M)$ of $g$:
 $$
  q(\,R\,)\;\;=\;\;\tfrac{\kappa}{8}\ .
 $$
 This assertion is essentially equivalent to the well--known
 Lichnerowicz--Schr\"odinger formula $D^2\,=\,\nabla^*\nabla\,+\,
 \frac{\kappa}{4}\,$ for the square of the Dirac operator $D$ acting
 on sections of the so called spinor bundle $\Sigma M$ (cf.~\cite{crelle}).
 In consequence $q(R)\,=\,\frac\kappa8$ for every geometric vector
 bundle occurring in the spinor bundle with respect to some holonomy
 reduction $\Hol\,M\,\subset\,\SO\,M$.

 \medskip\noindent
 In applications it is certainly useful to know several different
 presentations of the curvature endomorphism. A particularly elegant
 definition of $q(R)$ can be given by using the concept of conformal
 weights arising in conformal geometry. Consider for this purpose the
 conformal weight operator $B\,\in\,\Hom_\Hol(\,T\otimes T,\,\End V\,)$
 defined for all $X,\,Y\,\in\,T$ and $v\,\in\,V$ by:
 $$
  B_{X\,\otimes\,Y}\,v
  \;\;:=\;\;
  \pr_\hol(\,X\,\wedge\,Y\,)\,\star\,v\ ,
 $$
 Using the conformal weight operator we may write the curvature
 endomorphism in the form
 $$
  q(\,R\,)\,v
  \;\;=\;\;
  B(\,\nabla^2\,v\,)
  \;\;=\;\;
  \sum_{\mu\nu}B_{E_\mu\,\otimes\,E_\mu}\,\nabla^2_{E_\mu,\,E_\nu}v\ .
 $$
 for every section $v\,\in\,\Gamma(VM)$ of $VM$. Alternatively the conformal
 weight operator can be interpreted as a $\Hol$--equivariant endomorphism of
 the representation $T^*\otimes V$, as such it can be written in the form
 $B\,=\,\sum_\e b_\e\,\pr_\e$, where $\pr_\e$ are the projections to the
 isotypical components $V_\e\,\subset\,T^*\otimes V$. The eigenvalues or
 conformal weights $b_\e$ of $B$ can be computed easily by a very simple
 formula (cf.~\cite{compositio}). Via the relation $\pr_\e(\,\nabla^2\psi\,)
 \,=\,-P^*_\e P^{\phantom{*}}_\e\psi$ between the projections $\pr_\e$ and
 the hermitian squares of the generalized gradients $P_\e$ this calculation
 leads to the following {\it universal Weitzenb\"ock formula} studied
 extensively in \cite{compositio}:
 \begin{equation}\label{universal}
  q(R)
  \;\;=\;\;
  -\,\sum\,b_\e\,P^*_\e P^{\phantom{*}}_\e\ .
 \end{equation}
 For the generic Riemannian holonomy group $\SO\,T$ the conformal weight
 operator and the universal Weitzenb\"ock formula were first considered
 by P.~Gauduchon in \cite{gauduchon}, other uses of the conformal weight
 operator $B$ besides its use in conformal geometry can be seen
 in \cite{cgh} and \cite{homma}.
\subsection{Integral Representation}
\label{ir}
 Yet another presentation of the curvature endomorphism particularly
 useful under additional assumptions on the sectional curvatures of
 the Riemannian manifold $(\,M,\,g\,)$ writes $q(R)$ at a point $x\,\in\,M$
 as an integral over the Grassmannian of $2$--planes of $T_xM$ with respect
 to the Fubini--Study volume density. For a given cut off parameter
 $\L\,\in\,\R$ this integral representation of the curvature endomorphism
 is readily established by using the integration techniques for the
 sectional curvature discussed in \cite{secs} and reads
 \begin{eqnarray*}
  \lefteqn{\frac1{\mathrm{Vol}(\,\mathrm{Gr}_2T\,)}\,\genfrac(){0pt}{}{m+2}4
  \int_{\mathrm{Gr}_2T_xM}\Big(\,g(R_{X,Y}Y,X)\,-\,\L\,\Big)
  \,\pr_\hol(X\wedge Y)\star\pr_\hol(X\wedge Y)\star\,|\mathrm{vol}|}
  &&
  \\[4pt]
  &=&
  \Big(\frac\L{12}(m+2)(m+1)\,-\,\frac\kappa{12}\frac{m+4}m\Big)
  \,\Cas^\hol_{\L^2}\,+\,
  \sum_{\mu\nu}\pr_\hol(\Ric^\circ E_\mu\wedge E_\nu)\star
  \pr_\hol(E_\mu\wedge E_\nu)\star
  \\
  &&
  -\,\frac12\,q(\,R\,)\,+\,\frac1{48}\sum_{\lambda\mu\nu\rho}
  g(T\wedge T)(E_\lambda,E_\mu,E_\nu,E_\rho)\,\pr_\hol(E_\lambda\wedge E_\mu)
  \star\pr_\hol(E_\nu\wedge E_\rho)\star\ ,
 \end{eqnarray*}
 where $\kappa$ denotes the scalar curvature, $\Ric^\circ\,:=\,\Ric\,-\,
 \frac\kappa m\,\id$ the trace free part of the Ricci endomorphism and
 $\Cas^\hol_{\L^2}$ the Casimir operator of the Lie algebra bundle
 $\hol\,M$
 $$
  \Cas^\hol_{\L^2}
  \;\;:=\;\;
  -\,\tfrac12\,\sum_{\mu\nu}\pr_\hol (\,E_\mu\,\wedge\,E_\nu\,)\,\star\,
  \pr_\hol(\,E_\mu\,\wedge\,E_\nu\,)\,\star
 $$
 in $\L^2$--normalization. With respect to a $\Hol$--invariant scalar
 product or hermitian form on the representation $V$ corresponding to a
 geometric vector bundle $VM$ the endomorphism $\pr_\hol(X\wedge Y)\star$
 is skew symmetric or skew hermitian respectively with non--positive square
 in the sense of operators. Every upper bound $\mathrm{sec}\,\leq\,\Lambda$
 on the sectional curvatures of the Riemannian manifold $M$ thus leads in
 the torsion free case to a pointwise upper bound
 $$
  \frac12\,q(R)\;\;\leq\;\;
  \Big(\frac\L{12}(m+2)(m+1)\,-\,\frac\kappa{12}\frac{m+4}m\Big)
  \,\Cas^\hol_{\L^2}\,+\,
  \sum_{\mu\nu}\pr_\hol(\Ric^\circ E_\mu\wedge E_\nu)\star
  \pr_\hol(E_\mu\wedge E_\nu)\star
 $$
 for the curvature endomorphism $q(R)$, similarly every lower bound
 $\mathrm{sec}\,\geq\,\Lambda$ implies
 $$
  \frac12\,q(R)\;\;\geq\;\;
  \Big(\frac\L{12}(m+2)(m+1)\,-\,\frac\kappa{12}\frac{m+4}m\Big)
  \,\Cas^\hol_{\L^2}\,+\,
  \sum_{\mu\nu}\pr_\hol(\Ric^\circ E_\mu\wedge E_\nu)\star
  \pr_\hol(E_\mu\wedge E_\nu)\star
 $$
 in the sense of operators acting on the fiber $V_xM$ of a geometric vector
 bundles over a point $x\,\in\,M$. In \cite{aku}  this argument was used in
 the special case of geometric vector bundles $\S^p_\circ T^*M$ of
 trace--free symmetric $p$--tensors on a Riemannian manifold $M$
 with sectional curvature $\mathrm{sec} \le 0$ to show that the curvature
 endomorphism $q(R)$ is non--positive.
\subsection{The Standard Laplacian}
 The curvature term $q(R)$ is clearly functorial in the sense that it commutes
 with all morphism between geometric vector bundles. In order to obtain a
 second order differential operator enjoying the same functoriality property
 we simply add the rough Laplacian $\nabla^*\nabla$ and obtain the standard
 Laplace operator on geometric vector bundles:

 \begin{Definition}[Standard Laplace Operator]
 \hfill\label{slo}\break
  The standard Laplace operator acting on sections of a geometric vector
  bundle $VM$ is defined as the sum $\Delta\,=\,\nabla^*\nabla\,+\,q(R)$
  of the rough Laplacian and the curvature endomorphism $q(R)$.
 \end{Definition}

 \noindent
 The functorial nature of the standard Laplace operator for geometric vector
 bundles explains in a sense the work of many an author in differential
 geometry working for example in analogues of the Hodge decomposition of
 differential forms on K\"ahler manifolds: The standard Laplace operator
 commutes by definition with all morphisms of geometric vector bundles.
 Similar Laplace type operators with no or at most a limited functoriality
 have been present in the literature for a long time say in \cite{besse},
 section 1.~I, and in \cite{lichnerowicz} in the special case of tensor
 fields on Riemannian manifolds.

 \medskip
 {\bf Example 1:}
 In the Riemannian case $\Hol\,=\,\SO\;T$ and the geometric vector bundle
 $\L^pT^*M$ of differential forms of degree $p$ corresponding to the
 representation $V\,=\,\L^pT^*$ the standard Laplace operator $\Delta$
 coincides with the Hodge Laplace operator $\Delta_\Hodge\,:=\,
 d^*d\,+\,d d^*$. In fact the definition just reflects the classical
 Weitzenb\"ock formula for $\Delta_\Hodge$. However the functorial point
 of view is a decisive advantage in this example, because $\left.
 \Delta_\Hodge\right|_{VM}\,=\,\Delta$ for every geometric subbundle
 of the bundle of differential forms (cf.~\cite{crelle}) immediately
 implies the generalized Hodge decomposition of de Rham cohomology
 $$
  H^\bullet_{\mathrm{dR}}(\;M\;)
  \;\;=\;\;
  \bigoplus_V\Hom_\Hol(\,V,\,\L^\bullet T^*\,)
  \,\otimes\,\ker\,\Delta_V\ ,
 $$
 under arbitrary holonomy reductions $\Hol\,M\,\subset\,\SO\,M$, analogous
 decompositions hold true for every eigenvalue of the Hodge--Laplace operator
 $\Delta_\Hodge$ on the bundle of differential forms. The important point in
 this decomposition is that the restriction of $\Delta_\Hodge$ to a parallel
 subbundle $VM$ only depends on the representation $V$ and not on its
 embedding $V\,\subset\,\L^pT^*$. We note that a similar decomposition of 
 the de Rham cohomology is discussed in the book of D. Joyce
 (cf. Theorem 3.5.3 in \cite{joyce}).

 \medskip
 {\bf Example 2:} In the generic holonomy case $\Hol\,=\,\SO\,T$ the standard
 Laplace operator $\Delta$ coincides with the Lichnerowicz Laplacian
 $\Delta_L$  (cf.~\cite{besse}, 1.143, \cite{lichnerowicz}) on the
 geometric vector bundle $VM\,=\,\bigotimes^pT^*M$ of $p$--tensors.
 Especially interesting is the case of  trace--free symmetric $2$--tensors
 with associated representation $V\,=\,\S^2_0T^*$. In particular the space
 of infinitesimal Einstein deformations of an Einstein metric $g$ with
 Ricci curvature $\Ric\,=\,\frac\kappa m\,g$ and scalar curvature $\kappa$
 can be identified with the space of symmetric, trace and divergence free
 endomorphisms $H$ of $TM$ satisfying the eigenvalue equation $\Delta_LH
 \,=\,\frac{2\kappa}m\,H$ (cf.~\cite{besse}).

 \medskip
 {\bf Example 3:} On a Riemannian symmetric space $M\,=\,G/K$ with Riemannian
 metric $g$ induced by an invariant scalar product $g^\ext$ on the Lie algebra
 $\g$ geometric vector bundles are homogeneous and vice versa, moreover the
 standard curvature endomorphism $q(R)$ agrees with the Casimir operator of
 the isotropy algebra $\mathfrak{k}$. An easy calculation based on this fact
 mentioned above shows that the standard Laplace operator $\Delta$ is
 actually the Casimir operator of the isometry group $G$ on sections of
 a homogeneous vector bundle (cf.~Lemma 5.2 in \cite{au10}).

 \medskip\noindent
 Concerning the presentation of the curvature endomorphism $q(R)$ as a
 sum of hermitian squares of generalized gradients we remark that the
 universal Weitzenb\"ock formula \eqref{universal} extends directly to
 the standard Laplace operator, writing it as a linear combination of
 hermitian squares of generalized gradients with coefficients determined
 by the conformal weights $b_\e$
 \begin{equation}\label{delta}
  \Delta
   \;\;=\;\;
  \sum_\e(\,1\,-\,b_\e\,)\,P^*_\e P^{\phantom{*}}_\e\ ,
 \end{equation}
 which are as we have said before very easy to compute. This way to write
 the standard Laplace operators has been used in \cite{compositio} together
 with zero curvature term Weitzenb\"ock formulas to characterize for example
 all harmonic forms on $\G_2$ and $\Spin(7)$--manifolds.
\section{The Commutator Formula}
 \noindent
 In this section we will calculate the commutator of the standard Laplacian
 and a generalized gradient $D_\op$ from sections of a geometric vector bundle
 $VM$ to sections of a geometric vector bundle $\tilde VM$ over $M$. In order
 to compute the commutator of $\nabla^*\nabla$ and $D_\op$ it is convenient
 to recall the concept of iterated covariant derivatives for sections of a
 general vector bundle $VM$ over $M$ endowed with a connection $\nabla$ with
 respect to an auxiliary, not necessarily torsion free connection $\bar\nabla$
 on the tangent bundle $TM$. Specifically the second iterated covariant
 derivative is defined for every section $\psi\,\in\,\Gamma(\,VM\,)$
 by $\nabla^2_{X,\,Y}\psi\,:=\,\nabla_X(\nabla_Y\psi)\,-\,\nabla_{\bar
 \nabla_XY}\psi$ and much in the same spirit the third iterated covariant
 derivative reads:
 \begin{eqnarray*}
  \nabla^3_{X,\,Y,\,Z}\psi
  &:=&
  \nabla^{\phantom{1}}_X\nabla^2_{Y,\,Z}\psi
  \;-\;\nabla^2_{\bar\nabla_XY,\,Z}\psi
  \;-\;\nabla^2_{Y,\,\bar\nabla_XZ}\psi
  \\
  &=&
  \Big(\,\nabla_X\nabla_Y\nabla_Z
  -\nabla_X\nabla_{\bar\nabla_YZ}
  -\nabla_{\bar\nabla_XY}\nabla_Z
  +\nabla_{\bar\nabla_{\bar\nabla_XY}Z}
  -\nabla_Y\nabla_{\bar\nabla_X Z}
  +\nabla_{\bar\nabla_Y\bar\nabla_XZ}\,\Big)\psi\ .
 \end{eqnarray*}
 Our calculation of the commutator $[\,\nabla^*\nabla,\,D_\op\,]$
 relies on the following two identities for third iterated covariant
 derivatives known collectively as the {\it Ricci identities}
 (cf.~\cite{besse}, Corollary 1.22):
 \begin{eqnarray*}
  \nabla^3_{X,Y,Z}\;-\;\nabla^3_{Y,X,Z}
  &=&
  R_{X,\,Y}\nabla_Z
  \;-\;\nabla_{R_{X,Y}Z}
  \;-\;\nabla^2_{T(X,Y),Z}\ ,
  \\[4pt]
  \nabla^3_{Y,X,Z}\;-\;\nabla^3_{Y,Z,X}
  &=&
  (\nabla_Y R)_{X,Z}\;+\;R_{X,Z}\nabla_Y
  \;-\;\nabla_{(\bar\nabla_YT)(X,Z)}
  \;-\;\nabla^2_{Y,T(X,Z)}\ .
 \end{eqnarray*}

 \noindent
 Adding these two Ricci identities together while setting $Z\,=\,Y$ results
 in the identity:
 \begin{equation}\label{key}
  \begin{array}{lcl}
   \hbox to0pt{$\nabla^3_{X,Y,Y}\;-\;\nabla^3_{Y,Y,X}$\hss}\qquad
   &&
   \\[8pt]
   &=&
   2\,R_{X,Y}\,\nabla_Y
   -\nabla_{R_{X,Y}Y}
   -(\nabla_YR)_{Y,X}
   -\nabla^2_{T(X,Y),Y}
   -\nabla^2_{Y,T(X,Y)}
   -\nabla_{(\bar\nabla_YT)(X,Y)}\ .
  \end{array}
 \end{equation}
 This is the key identity in the proof of the commutator formula for the
 standard Laplacian:

 \begin{Theorem}[Commutator Formula]
 \hfill\label{theorem}\break
  Let $(\,M,\,g\,)$ be a Riemannian manifold and let $\bar\nabla$ be a
  metric connection with parallel skew--symmetric torsion $T$ in the sense
  that the expression $\theta(X,Y,Z)\,:=\,g(T(X,Y),Z)$ defines a parallel
  differential form $\theta\,\in\,\Gamma(\,\L^3T^*M\,)$. Consider a
  generalized gradient
  $$
   D_\op:\;\;
   \Gamma(\,VM\,)
   \;\stackrel\nabla\longrightarrow\;\Gamma(\,T^*M\otimes VM\,)
   \;\stackrel\cong\longrightarrow\;\Gamma(\,TM\otimes VM\,)
   \;\stackrel{\sigma_\op}\longrightarrow\;\Gamma(\,\tilde VM\,)
  $$
  from a geometric vector bundle $VM$ to a geometric vector bundle $\tilde VM$
  associated to a $\Hol$--equivariant bilinear operation $\op:\,T\otimes V
  \longrightarrow\tilde V$ between the associated representations. The
  commutator of the standard Laplace operator $\Delta$ acting on both
  $VM$ and $\tilde VM$ with $D_\op$ reads:
  \begin{eqnarray*}
   [\;\Delta,\;D_\op\;]
   \;\;=\;\;
   \Delta_{\tilde V}\,\circ\,D_\op\;-\;D_\op\,\circ\,\Delta_V
   &=&
   -\,\sigma_\op(\;q(\,\bar\nabla R\,)\;-\;\delta R\,\star\;)
   \\
   &:=&
   -\,\sum_\lambda E_\lambda\,\op\,(\,q(\bar\nabla_{E_\lambda}R)
   \;-\;(\delta R)_{E_\lambda}\,\star\,)\ .
  \end{eqnarray*}
 \end{Theorem}
 
 \proof
 In terms of the parallel extension $\op:\,TM\otimes VM
 \longrightarrow\tilde VM$ of $\op$ we may write $D_\op$ as
 $$
  D_\op
  \;\;:=\;\;
  \sigma_\op\,\circ\,\nabla
  \;\;=\;\;
  \sum_\lambda E_\lambda\,\op\,\nabla_{E_\lambda}\ .
 $$
 with a local orthonormal basis $\{\,E_\lambda\,\}$. Using the parallelity
 of $\op$ we find for all $X\,\in\,\Gamma(\,TM\,)$
 \begin{eqnarray*}
  \nabla_X\,\circ\,D_\op
  &=&
  \sum_\lambda\Big(\,(\bar\nabla_XE_\lambda)\,\op\,\nabla_{E_\lambda}
  \;+\;E_\lambda\,\op\,\nabla_X\nabla_{E_\lambda}\;\Big)
  \\
  &=&
  \sum_\lambda E_\lambda\,\op\,\Big(\;-\;\nabla_{\bar\nabla_XE_\lambda}
  \;+\;\nabla_X\nabla_{E_\lambda}\;\Big)
  \;\;=\;\;
  \sum_\mu E_\lambda\,\op\,\nabla^2_{X,\,E_\lambda}\ ,
 \end{eqnarray*}
 where in the second line we used the parallelity $\sum_\lambda(\bar\nabla_X
 E_\lambda)\otimes E_\lambda\,=\,-\sum_\lambda E_\lambda\otimes(\bar\nabla_X
 E_\lambda)$ of the cometric tensor with respect to a metric connection.
 Applying this identity twice we find:
 \begin{eqnarray*}
  \nabla^*\nabla\,\circ\,D_\op
  &=&
  -\,\sum_{\lambda\mu}E_\lambda\,\op\,\nabla^3_{E_\mu,\,E_\mu,\,E_\lambda}
  \\
  D_\op\,\circ\,\nabla^*\nabla
  &=&
  -\,\sum_{\lambda\mu}E_\lambda\,\op\,\nabla^3_{E_\lambda,\,E_\mu,\,E_\mu}
 \end{eqnarray*}
 In consequence the commutator $[\,\nabla^*\nabla,\,D_\op\,]$ equals
 the sum of the key identity (\ref{key})
 \begin{eqnarray*}
  [\,\nabla^*\nabla,\,D_\op\,]
  &=&
  \sum_{\lambda\mu}E_\lambda\,\op\,
  \Big(\;\nabla^3_{E_\lambda,\,E_\mu,\,E_\mu}
  \;-\;\nabla^3_{E_\mu,\,E_\mu,\,E_\lambda}\;\Big)
  \\
  &=&
  2\,\sum_{\lambda\mu}E_\lambda\,\op\,R_{E_\lambda,\,E_\mu}\,\nabla_{E_\mu}
  \,-\,
  \sum_{\lambda\mu}E_\lambda\,\op\,\nabla_{R_{E_\lambda,\,E_\mu}E_\mu}
  \,-\,
  \sum_{\lambda\mu}E_\lambda\,\op\,
  (\,\bar\nabla_{E_\mu}R\,)_{E_\mu,\,E_\lambda}\ .
 \end{eqnarray*}
 over a local orthonormal basis, because the torsion $\bar\nabla T\,=\,0$
 is parallel by assumption and the other two terms involving $T$ cancel out
 due to $\sum_\mu(\,T(X,E_\mu)\otimes E_\mu\,+\,E_\mu\otimes T(X,E_\mu)\,)
 \,=\,0$ for all $X$, after all $Y\longmapsto T(X,Y)$ is a skew symmetric
 endomorphism of $TM$ for skew--symmetric torsion $T$. In terms of the Ricci
 endomorphism $\Ric\,X\,:=\,\sum_\mu R_{X,\,E_\mu}E_\mu$ and the divergence
 $(\delta R)_X\,:=\,-\sum_\mu(\nabla_{E_\mu}R)_{E_\mu,\,X}$ of the curvature
 tensor this result can be written:
 \begin{equation}\label{eins}
  [\,\nabla^*\nabla,\,D_\op\,]
  \;\;=\;\;
  2\,\sum_{\lambda\mu}E_\lambda\,\op\,R_{E_\lambda,\,E_\mu}\,\nabla_{E_\mu}
  \;-\;\sum_\lambda E_\lambda\,\op\,\nabla_{\Ric\,E_\lambda}
  \;+\;\sum_\lambda E_\lambda\,\op\,(\,\delta R\,)_{E_\lambda}\star\ .
 \end{equation}
 In a second step we calculate the commutator of the curvature endomorphism
 $q(R)$ characteristic for the standard Laplace operator $\Delta$ with
 the generalized gradient $D_\op$. Recalling the alternative definition
 of the curvature term $q(R)\,=\,\sum_\alpha\X_\alpha\star R(\,\X_\alpha\,)
 \,\star$ as a sum over a local orthonormal basis $\{\,X_\alpha\,\}$ of the
 holonomy algebra bundle $\hol\,M\,\subset\,\L^2TM$ we find
 \begin{eqnarray*}
  [\;q(\,R\,),\;D_\op\;]
  &=&
  \sum_\lambda\Big(\,q(\,R\,)\,E_\lambda\,\Big)\,\op\,\nabla_{E_\lambda}
  \;+\;\sum_\lambda E_\lambda\,\op\,\Big(\,q(\,R\,)\,\nabla_{E_\lambda}\,
  \Big)
  \\
  &&
  \quad+\;2\,\sum_{\lambda\mu\alpha}g(\X_\alpha\star E_\lambda,E_\mu)
  \,E_\mu\,\op\,R(\,\X_\alpha\,)\,\star\,\nabla_{E_\lambda}
  \;-\;\sum_\lambda E_\lambda\,\op\,\nabla_{E_\lambda}\,q(\,R\,)
  \\[3pt]
  &=&
  \sum_\lambda\Ric\,E_\lambda\,\op\,\nabla_{E_\lambda}
  \;-\;\sum_\lambda E_\lambda\,\op\,q(\,\bar\nabla_{E_\lambda}R\,)
  \\
  &&
  \quad+\;2\,\sum_{\lambda\mu}E_\mu\,\op\,R\,\Big(\,\sum_\alpha
  g(\X_\alpha,E_\lambda\wedge E_\mu)\,\X_\alpha\,\Big)
  \,\star\,\nabla_{E_\lambda}
  \\
  &=&
  \sum_\lambda\Ric\,E_\lambda\,\op\,\nabla_{E_\lambda}
  \;-\;\sum_\lambda E_\lambda\,\op\,q(\,\bar\nabla_{E_\lambda}R\,)
  \;+\;2\,\sum_{\lambda\nu}E_\nu\,\op\,R_{E_\lambda,\,E_\nu}\,
  \nabla_{E_\lambda}\ ,
 \end{eqnarray*}
 because the curvature endomorphism $q(R)$ equals the Ricci endomorphism
 on the tangent bundle and $R(\,\pr_\hol(E_\lambda\wedge E_\mu)\,)\star
 \,=\,R(\,E_\lambda\wedge E_\mu\,)\star\,=\,R_{E_\lambda,\,E_\mu}$ as the
 curvature endomorphism factorizes over the orthogonal projection to the
 holonomy algebra bundle $\hol\,M\,\subset\,\L^2TM$. Adding the commutator
 $[\,\nabla^*\nabla,\,D_\op\,]$ from equation (\ref{eins}) to the commutator
 $[\,q(R),\,D_\op\,]$ calculated above and using the symmetry $\sum_\lambda
 (\Ric\,E_\lambda)\otimes E_\lambda\,=\,\sum_\lambda E_\lambda\otimes
 (\Ric\,E_\lambda)$ of the Ricci curvature endomorphism we eventually
 obtain the commutator formula:
 $$
  [\;\Delta,\;D_\op\;]
  \;\;=\;\;
  -\;\sum_\lambda E_\lambda\,\op\,q(\,\bar\nabla_{E_\lambda}R\,)
  \;+\;\sum_\lambda E_\lambda\,\op\,(\,\delta R\,)_{E_\lambda}\,\star\ .
 $$
 \vskip-28pt
 \qed

 \begin{Corollary}[Commutator Formula for First Order Operators]
 \hfill\label{cf1o}\break
  Every natural first order differential operator $D:\,\Gamma(\,VM\,)
  \longrightarrow\Gamma(\,\tilde VM\,)$ between geometric vector bundles
  can be written as a sum $D\,=\,D_\op\,+\,F$ of the generalized gradient
  associated to its principal symbol $\sigma_\op:\,T^*M\otimes VM
  \longrightarrow\tilde VM$ and a homomorphism $F:\,VM\longrightarrow\tilde
  VM$ of geometric vector bundles. In consequence we find
  $[\,\Delta,\,D\,]\,=\,[\,\Delta,\,D_\op\,]$.
 \end{Corollary}

 \noindent
 Because both error terms $\sum_\lambda E_\lambda\op
 q(\bar\nabla_{E_\lambda}R)$ and $\sum_\lambda E_\lambda\op
 (\delta R)_{E_\lambda}\star$ on the right hand side of the commutator
 formula of Theorem \ref{theorem} are the images of the covariant
 derivative $\bar\nabla R$ of the curvature tensor under homomorphisms
 of geometric vector bundles, we conclude:

 \begin{Corollary}[Simple Vanishing Criterion]
 \hfill\label{corollary}\break
  Let $\Hol$ be the holonomy group of the metric connection $\bar\nabla$ and
  suppose that $\bar\nabla R$ is a section of a geometric vector bundle
  $CM$ with corresponding $\Hol$--representation $C$. A sufficient condition
  for the vanishing of the commutator $[\,\Delta,\,D\,]\,=\,0$ of the standard
  Laplace operator with a natural first order differential operator
  $D$ from a geometric vector bundle $VM$ to a geometric vector bundle
  $\tilde VM$ is the vanishing of the space of $\Hol$--equivariant
  homomorphisms:
  $$
   \Hom_\Hol(\;C,\;\Hom(\,V,\,\tilde V\,)\;)
   \;\;=\;\;
   \{\,0\,\}\ .
  $$
 \end{Corollary}

 \noindent
 For later uses it may be helpful to recall the well--known fact that for
 the Levi--Civita connection $\bar\nabla$ for the Riemannian metric $g$
 the divergence $\delta R$ of the Riemannian curvature tensor $R$ relates
 via $\delta R\,=\,d^{\bar\nabla}\Ric$ to the covariant derivative of the
 Ricci curvature. In consequence $\delta R$ vanishes on all manifolds
 with parallel Ricci tensor and the error term in the commutator formula
 involves only $q(\bar\nabla R)$. According to an observation of Gray
 \cite{gray} manifolds with parallel Ricci tensor are locally the Cartesian
 product of Einstein manifolds and conversely every local product of Einstein
 manifolds has parallel Ricci curvature provided all $2$--dimensional
 Einstein factors have constant scalar curvature.
\section{Examples}
\subsection{Differential and Codifferential}
 Consider a Riemannian manifold $(\,M,\,g\,)$ and the Hodge--Laplace
 operator $\Delta_\Hodge\,:=\,d^*d\,+\,dd^*$ acting on the space of
 $p$--forms $\Gamma(\,\L^pT^*M\,)$. Because both $d$ and $d^*$ are
 boundary operators $d^2\,=\,0\,=\,d^{*2}$ the Hodge--Laplace operator
 commutes with both $d$ and $d^*$ due to the rather trivial argument:
 $$
  [\;\Delta_\Hodge,\;d\;]
  \;\;=\;\;
  d^*d^2\,+\,dd^*d\,-\,dd^*d\,-\,d^2d^*
  \;\;=\;\;
  0\ .
 $$
 The commutator formula of Theorem \ref{theorem} can be interpreted as a
 vast generalization of this simple observation forming the philosophical
 underpinning of Hodge theory. Taking $\bar\nabla$ to the Levi--Civita
 connection on $TM$ and $\nabla$ to be the induced connection on the
 geometric vector bundle $\L^pT^*M$ of $p$--forms we may write the exterior
 derivative $d\,=\,\sigma_\wedge\circ\nabla$ as a generalized gradient
 corresponding to the wedge product $\sigma_\wedge:\,T\otimes\L^pT^*
 \longrightarrow\L^{p+1}T^*,\,X\otimes\omega\longmapsto X^\flat\wedge\omega$.
 The commutator $[\,\Delta_\Hodge,\,d\,]$ thus falls under the ambit of
 Theorem \ref{theorem}, because the standard Laplace operator $\Delta\,=\,
 \Delta_\Hodge$ agrees with the Hodge--Laplace operator in this setup as
 we remarked previously. In order to determine this commutator we verify
 first of all the equation:
 \begin{eqnarray*}
  \sigma_\wedge(q(\bar\nabla R))
  &=&
  -\sum_{\lambda\mu\nu}
  E_\lambda^\flat\wedge\,E_\mu^\flat\wedge\,E_\nu\,\ins\,
  (\nabla_{E_\lambda}R)_{E_\mu,\,E_\nu}
  \\
  &=&
  -\sum_{\lambda\mu\nu}
  E_\nu\,\ins\,E_\lambda^\flat\wedge\,E_\mu^\flat\wedge\,
  (\nabla_{E_\lambda}R)_{E_\mu,\,E_\nu}
  +\sum_{\lambda\mu}
  E_\mu^\flat\wedge\,(\nabla_{E_\lambda}R)_{E_\mu,\,E_\lambda}
  \;\,=\;\,
  \sigma_\wedge(\delta R\star)\ .
 \end{eqnarray*}
 The omission of the first sum is justified by using the second Bianchi
 identity in the calculation
 \begin{eqnarray*}
  \sum_{\lambda\mu}
  E_\lambda^\flat\wedge\,E^\flat_\mu\wedge\,
  (\,\nabla_{E_\lambda}R\,)_{E_\mu,\,X}
  &=&
  +\,\tfrac12\,\sum_{\lambda\mu}E_\lambda^\flat\wedge\,E_\mu^\flat\wedge\,
  \left(\;(\,\nabla_{E_\lambda}R\,)_{E_\mu,\,X}\;-\;
  (\,\nabla_{E_\mu}R\,)_{E_\lambda,\,X}\;\right)
  \\
  &=&
  -\,\tfrac12\,\sum_{\lambda\mu}
  E_\lambda^\flat\wedge\,E_\mu^\flat\wedge\,
  (\,\nabla_XR\,)_{E_\lambda,\,E_\mu}
  \;\;=\;\;
  0
 \end{eqnarray*}
 for all $X\,\in\,\Gamma(\,TM\,)$, where $\sum_{\lambda\mu}E^\flat_\lambda
 \wedge E^\flat_\mu\wedge R'_{E_\lambda,\,E_\mu}\,=\,0$ holds true for
 every algebraic curvature tensor and thus in particular for $R'\,=\,
 \nabla_XR$ due to the first Bianchi identity. In passing we remark that
 this consequence of the first Bianchi identity arises from the proof of
 the classical Weitzenb\"ock formula for the Hodge--Laplace operator
 $\Delta_\Hodge$. The error term $\sigma_\wedge(\,q(\nabla R)\,-\,
 \delta R\star\,)\,=\,0$ in the commutator formula of Theorem \ref{theorem}
 thus vanishes for the exterior derivative and $\Delta_\Hodge\,=\,\Delta$
 commutes with $d$. Mutatis mutandis this argument establishes the vanishing
 of the commutator $[\,\Delta_\Hodge,\,d^*\,]\,=\,0$ as well.
\subsection{Symmetric and Normal Homogeneous Spaces}
 Symmetric spaces are essentially Riemannian manifolds $(\,M,\,g\,)$ with
 parallel Riemannian curvature tensor $R$ with respect to the Levi--Civita
 connection $\nabla$. The error terms in the commutator formula of Theorem
 \ref{theorem} corresponding to $q(\nabla R)$ and $\delta R$ thus
 both vanish, in consequence the standard Laplacian $\Delta$ commutes
 with all equivariant first order differential operators between
 geometric, i.~e.~homogeneous vector bundles on $M$ vindicating the
 identification of $\Delta$ with the Casimir operator of the isometry
 group of a symmetric space $M$.

 Interestingly the same statement is true on compact normal homogeneous
 spaces $(\,M,\,g\,)$ characterized by a Riemannian metric induced from
 an invariant scalar product on the Lie algebra of the isometry group. The
 metric connection of choice under this assumption is not the Levi--Civita
 connection $\nabla$ however, but the reductive connection $\bar\nabla$
 characterized by having parallel skew symmetric torsion and parallel
 curvature tensor $R$. With this proviso the standard Laplace operator
 on geometric vector bundles can be identified with the Casimir operator
 of the isometry group of a normal homogeneous space $M$ (cf.~Lemma 5.2
 in \cite{au10}).
\subsection{The Rarita--Schwinger Operator}
 Let $(\,M,\,g\,)$ be a Riemannian spin manifold with the spin structure
 defined by a lift $\Spin\,M$ of the principal $\SO\,T$--bundle $\SO\,M$
 of oriented orthonormal frames to the structure group $\Spin\,T$ as before.
 Because the Clifford multiplication $\bullet:\,T\otimes\Sigma\longrightarrow
 \Sigma$ for the spinor representation $\Sigma$ is equivariant under
 $\Spin\,T$, its kernel $\Sigma_{\frac32}\,\subset\,T\otimes\Sigma$ is
 a $\Spin\,T$--representation with associated geometric vector bundle
 $\Sigma_{\frac32}M$ on $M$. In passing we remark that $\Sigma_{\frac32}$
 equals the Cartan summand in the tensor product $T\otimes\Sigma$ in odd
 dimensions, in even dimensions $\Sigma_{\frac32}\,=\,\Sigma^+_{\frac32}
 \oplus\Sigma^-_{\frac32}$ decomposes like the spinor representation itself.
 The {\em Rarita--Schwinger operator} is defined as the generalized gradient
 $$
  D_\op:\qquad\Gamma(\,\Sigma_{\frac32}M\,)
  \;\longrightarrow\;\Gamma(\,\Sigma_{\frac32}M\,)
 $$
 associated to the $\Spin\,T$--equivariant homomorphism $\sigma_\op:\,
 T\otimes\Sigma_{\frac32}\longrightarrow\Sigma_{\frac32}$ given by Clifford
 multiplication in the $\Sigma$--factor followed by projection 
 $X\op\psi\,:=\,\pr_{\Sigma_{\frac32}}(\,(\id\otimes X\bullet)\psi\,)$.
 It is well--known, compare for example \cite{galaev} or \cite{salamon},
 that the covariant derivative of the curvature tensor of an Einstein
 manifold is a section $\nabla R\,\in\,\Gamma(\,CM\,)$ of the geometric
 vector bundle corresponding to the Cartan summand $C\,:=\,\S^{3,2}_\circ
 T^*\,\subset\,\S^3_\circ T^*\otimes\S^2_\circ T^*$ in the tensor product
 of harmonic cubic and quadratic polynomials on $T$; and it is easily
 verified that this Cartan summand does not occur in the endomorphisms
 of $\Sigma_{\frac32}$. In consequence $\Delta$ commutes $[\,\Delta,\,
 D_\op\,]\,=\,0$ with the Rarita--Schwinger operator $D_\op$ on Einstein
 spin manifolds due to Corollary \ref{corollary}. A direct proof of this
 commutator formula can be found in \cite{homma1}.
\subsection{Nearly K\"ahler manifolds}
 A six dimensional nearly K\"ahler manifold $(\,M,\,g,\,J\,)$ is a Riemannian
 manifold of dimension $6$ endowed with an orthogonal almost complex structure
 $J$ satisfying the integrability condition $(\nabla_X J)X\,=\,0$ for the
 Levi--Civita connection $\nabla$ and all vector fields $X\,\in\,\Gamma(TM)$.
 In this situation it is better to replace the Levi--Civita connection
 by the {\em canonical hermitian connection} $\bar\nabla_X\,:=\,\nabla_X
 \,+\,\frac12\,J(\nabla_XJ)$ on the tangent bundle, which is a metric
 connection with parallel, skew--symmetric torsion in the sense
 of Definition \ref{cpt} and makes $J$ parallel $\bar\nabla J\,=\,0$.
 In consequence the holonomy group of $\bar\nabla$ is a subgroup of
 $\Hol\,=\,\SU(T,J,\theta)\,\subset\,\SO\;T$ unless $M$ is actually
 K\"ahler with $\nabla J\,=\,0$.

 Considering only the case of strictly nearly K\"ahler manifolds with
 $\nabla J\,\neq\,0$ we conclude that the defining $3$--dimensional
 representation of $\SU(T,J,\theta)\,\cong\,\SU(3)$ gives rise to a
 geometric vector bundle $EM$ on $M$ such that $TM\otimes_\R\C\,=\,
 EM\oplus\bar EM$. Due to $\L^\bullet E\,\cong\,\C\oplus E\oplus\bar E
 \oplus\C$ etc.~the geometric vector bundle of complex--valued differential
 forms on $M$ decomposes
 $$
  \L^\bullet T^*M\otimes_\R\C
  \;\;\cong\;\;
  (\,\C M\,\oplus\,EM\,\oplus\,\bar EM\,\oplus\,\C M\,)
  \,\otimes\,(\,\C M\,\oplus\,\bar EM\,\oplus\,EM\,\oplus\,\C M\,)
 $$
 into copies of the trivial bundle $\C M$, $EM$ and $\bar EM$, the
 complexified holonomy algebra bundle $\su\,M\otimes_\R\C$ and $\S^2EM$
 as well as $\S^2\bar EM$. On the other hand it is well--known that the
 curvature tensor of the canonical hermitian connection $\bar\nabla$ can
 be written as a sum $R\,=\,\frac\kappa{40}\,R^\circ\,+\,R^{\mathrm{CY}}$
 of a parallel standard curvature tensor $R^\circ$ and a curvature tensor
 of Calabi--Yau type, i.~e.~a real section $R^{\mathrm{CY}}\,\in\,\Gamma
 (\,\S^2EM\otimes_\circ\S^2\bar EM\,)$ of the vector bundle $\S^2EM
 \otimes_\circ\S^2\bar EM$. Its covariant derivative then is a real
 section $\bar\nabla R\,\in\,\Gamma(CM)$ of the vector bundle corresponding
 to the representation $C\,=\,\S^3E\otimes_\circ\S^2\bar E\,\oplus\,\S^2E
 \otimes_\circ\S^3\bar E$, which is simply too complex to occur in the
 endomorphisms of $\L^\bullet T^*\otimes_\R\C$:

 \begin{Proposition}[Commutator Formulas for Nearly K\"ahler Manifolds]
 \hfill\label{cnk}\break
  Let $(\,M,\,g,\,J\,)$ be a six dimensional nearly K\"ahler manifold with
  canonical hermitian connection $\bar\nabla$. The standard Laplacian
  $\Delta\,=\,\bar\nabla^*\bar\nabla\,+\,q(\bar R)$ acting on differential
  forms commutes with the exterior derivative $d$ and the codifferential
  $d^*$ and thus with the Hodge--Laplace operator:
  $$
   [\;\Delta,\;d\;]
   \;\;=\;\;
   0
   \qquad\qquad
   [\;\Delta,\;d^*\;]
   \;\;=\;\;
   0
   \qquad\qquad
   [\;\Delta,\;\Delta_\Hodge\;]
   \;\;=\;\;
   0\ .
  $$
  In fact $\Delta$ commutes with all natural first order differential
  operators on differential forms.
 \end{Proposition}
\subsection{$\G_2$-- and $\Spin(7)$--Manifolds}
 Consider a Riemannian manifold $(\,M,\,g\,)$ of dimension seven with
 holonomy $\G_2$. The irreducible geometric vector bundles associated to
 this holonomy reduction to $\Hol\,=\,\G_2\,\subset\,\SO\,T$ correspond to
 the irreducible representations $V_{[\,a,\,b\,]}$ of $\G_2$, which are
 parametrized by their highest weight, a linear combination $a\omega_1\,+\,
 b\omega_2$ with integer coefficients $a,\,b\,\geq\,0$ of the two fundamental
 weights $\omega_1$ and $\omega_2$ corresponding respectively to the
 $7$--dimensional isotropy representation $T$ and the adjoint representation
 $\hol\,=\,\g_2$. It is possible to check that the Riemannian curvature tensor
 $R$ is a section of the geometric vector bundle $V_{[\,0,\,2\,]}M$ and that
 its covariant derivative $\nabla R$ is a section of $V_{[\,1,\,2\,]}M$,
 compare \cite{salamon}, page 162. The representation of $\G_2$ on the exterior
 algebra $\L^\bullet T^*$ of alternating forms on $T$ splits on the other
 hand into copies of the trivial representation $\R$, the isotropy
 representation $T$, the holonomy representation $\hol\,=\,\g_2$ and
 $\L^3_{27}\,:=\,V_{[\,2,\,0\,]}$. Considering the decompositions
 \begin{eqnarray*}
  \hbox to18pt{\hfill$T$\hfill}
  \otimes\hbox to18pt{\hfill$\g_2$\hfill}
  &=&
  V_{[\,1,\,0\,]} \oplus V_{[\,2,\,0\,]} \oplus V_{[\,1,\,1\,]}\ ,
  \\[3pt]
  \hbox to18pt{\hfill$T$\hfill}
  \otimes\hbox to18pt{\hfill$\L^3_{27}$\hfill}
  &=&
  V_{[\,1,\,0\,]}\oplus V_{[\,2,\,0\,]}\oplus V_{[\,3,\,0\,]}
  \oplus V_{[\,0,\,1\,]}\oplus V_{[\,1,\,1\,]}\ ,
  \\[3pt]
  \hbox to18pt{\hfill$\g_2$\hfill}
  \otimes\hbox to18pt{\hfill$\g_2$\hfill}
  &=&
  V_{[\,0,\,0\,]} \oplus V_{[\,2,\,0\,]} \oplus V_{[\,3,\,0\,]} \oplus
  V_{[\,0,\,1\,]} \oplus V_{[\,0,\,2\,]}\ ,
  \\[3pt]
  \hbox to18pt{\hfill$\g_2$\hfill}
  \otimes\hbox to18pt{\hfill$\L^3_{27}$\hfill}
  &=&
  V_{[\,1,\,0\,]} \oplus V_{[\,2,\,0\,]} \oplus V_{[\,3,\,0\,]} \oplus
  V_{[\,0,\,1\,]} \oplus V_{[\,1,\,1\,]} \oplus V_{[\,2,\,1\,]}\ ,
  \\[3pt]
  \hbox to18pt{\hfill$\L^3_{27}$\hfill}
  \otimes\hbox to18pt{\hfill$\L^3_{27}$\hfill}
  &=&
  V_{[\,0,\,0\,]} \oplus V_{[\,1,\,0\,]} \oplus 2\,V_{[\,2,\,0\,]} \oplus
  V_{[\,3,\,0\,]} \oplus V_{[\,4,\,0\,]} \oplus V_{[\,0,\,1\,]}\oplus
  2\,V_{[\,1,\,1\,]} \oplus V_{[\,2,\,1\,]} \oplus V_{[\,0,\,2\,]}
 \end{eqnarray*}
 we conclude that $\nabla R$ cannot result in a homomorphism between any two
 irreducible components of the exterior algebra $\L^\bullet T^*$ of alternating
 forms on $T$, in turn Corollary \ref{corollary} tells us that the standard
 Laplace operator $\Delta$ commutes with every natural first order
 differential operator on differential forms. In particular $\Delta$
 commutes $[\,\Delta,\,d_c\,]\,=\,0$ with the modified differential
 $d_c$ introduced by Verbitsky in \cite{verbitsky}. A very similar
 result holds true in the case of $\Spin(7)$--holonomy, the details
 of this argument are left to the reader.
\subsection{Quaternion--K\"ahler Manifolds}
 A Riemannian manifold $(\,M,\,g\,)$ of dimension $4n$ divisible by $4$
 with holonomy group contained in $\Sp(1)\cdot\Sp(n)\,\subset\,\SO(4n)$
 is called a quaternion--K\"ahler manifold. In order to describe the geometric
 vector bundles on a quaternion--K\"ahler manifold $M$ associated to this
 holonomy reduction we will denote the defining $2$-- and $2n$--dimensional
 complex representations of $\Sp(1)$ and $\Sp(n)$ respectively by $H$ and
 $E$. In general $H$ and $E$ do not give rise to geometric vector bundles
 on $M$, because neither representation extends to a representation of
 $\Hol\,=\,\Sp(1)\cdot\Sp(n)$, however all totally even powers of $H$
 and $E$ do. The complexified tangent bundle for example corresponds to
 the geometric vector bundle $TM\otimes_\R\C\,=\,HM\otimes EM$.

 The Riemannian curvature tensor $R$ of a quaternion--K\"ahler manifold $M$
 can be written as the sum of the curvature tensor of the quaternionic
 projective space $\H P^n$ with the same scalar curvature and a curvature
 tensor of hyperk\"ahler type $R\,=\,R^{\H P^n}\,+\,R^{\mathrm{hyper}}$,
 i.~e.~a real section $R^{\mathrm{hyper}}\,\in\,\Gamma(\,\S^4EM\,)$ of the
 geometric vector bundle $\S^4EM$. Working out the details of the second
 Bianchi identity we observe that the covariant derivative $\nabla R$ of
 the curvature tensor is a real section of the geometric vector bundle
 $HM\otimes\S^5 EM$ (cf.~\cite{crelle}). The geometric vector bundle
 $\L^\bullet T^*M\otimes_\R\C$ of complex--valued alternating forms on
 $M$ decomposes in a rather complicated way into a sum of geometric vector
 bundles of the form $\S^kHM\otimes\L^{a,b}_\circ EM$ with $k\,\geq\,0$
 and $n\,\geq\,a\,\geq\,b\,\geq\,0$ (cf.~\cite{gregor}), where $\L^{a,b}_0E
 \,\subset\,\L^a_\circ E\otimes\L^b_\circ E$ is the Cartan summand in the
 tensor product of the kernels $\L^a_\circ E\,\subset\,\L^aE$ and
 $\L^b_\circ E\,\subset\,\L^bE$ of the contraction with the symplectic
 form. Generically there are ten generalized gradients defined for
 sections of the geometric vector bundles $\S^kHM\otimes\L^{a,b}_\circ EM$
 and the standard Laplace operator $\Delta$ commutes with at least eight
 of these generalized gradients:

 \begin{Proposition}[Commutator Formula for Quaternion--K\"ahler Manifolds]
 \hfill\label{qkf}\break
  Every generalized gradient $D_\op:\,\Gamma(\,\S^kHM\otimes\L^{a,b}_\circ
  EM\,)\longrightarrow\Gamma(\,VM\,)$ to sections of a geome\-tric vector
  bundle $VM$ on a quaternion--K\"ahler manifold $M$ corresponding to a
  representation of the form $V\,=\,\S^{k\pm1}H\otimes\L^{a\pm 1,b}_\circ E$
  or $V\,=\,\S^{k\pm1}H\otimes\L^{a,b\pm1}_\circ E$ commutes with $\Delta$.
 \end{Proposition}

 \proof
 Because all quaternion--K\"ahler manifolds are automatically Einstein,
 we only need to consider the error term $\sigma_\op(\,q(\nabla R)\,)$ in
 the commutator formula of Theorem \ref{theorem}, where $\sigma_\op$ is the
 parallel vector bundle extension of the $\Sp(1)\cdot\Sp(n)$--equivariant
 isotypical projection corresponding to $D_\op$. Since $\sigma_\op(\,
 q(\nabla R\,)$ is a homomorphism of vector bundles parametrized by the
 section $\nabla R\,\in\,\Gamma(\,HM\otimes\S^5EM\,)$, we may replace all
 geometric vector bundles by their respective representations reducing the
 problem to the problem to show that $\sigma_\op(\,q(R')\,)\,=\,0$ for all
 $R'\,\in\,H\otimes\S^5E$ or equivalently for all linear generators $R'\,=\,
 h\otimes\frac1{5!}e^5$ with arbitrary $h\,\in\,H$, $e\,\in\,E$. In a similar
 vein $q(\frac1{4!}e^4)\,\sim\,e\wedge e^\flat\ins\otimes e\wedge e^\flat\ins$
 according to \cite{crelle} for every subrepresentation of $\L\,E\otimes\L\,E$
 and so in particular for $\L^{a,b}_\circ E$. In turn the homomorphism
 $$
  \sigma_\op(\,q(\,h\otimes\tfrac1{5!}e^5\,)\,):\qquad
  \S^kH\,\otimes\,\L^{a,b}_\circ E\;\longrightarrow\;\S^{k-1}H\,
  \otimes\,\L^{a+1,b}_\circ E
 $$
 for example equals the trivial homomorphism $\sigma_\op(\,q(h\otimes
 \frac1{5!}e^5)\,)\,=\,0$, because
 \begin{eqnarray*}
  \sigma_\op\Big(\,q(h\otimes\tfrac1{5!}e^5)\,\Big)\,
  (\,\alpha\otimes\eta\,)
  &:=&
  (\,h^\flat\,\ins\,\alpha\,)\otimes
  \pr_{\L^{a+1,b}_\circ E}\Big(\,
  (\,e\wedge\otimes\id\,)\,q(\tfrac1{4!}e^4)\,\eta\,\Big)
  \\
  &\sim&
  (\,h^\flat\,\ins\,\alpha\,)\otimes
  \pr_{\L^{a+1,b}_\circ E}\Big(\;(\,e\wedge\,e\wedge\,e^\flat\ins\,\otimes
  e\wedge e^\flat\ins\,)\;\eta\;\Big)
  \;\;=\;\;
  0
 \end{eqnarray*}
 for all $\alpha\,\in\,\S^kH$ and $\eta\,\in\,\L^{a,b}_\circ E$, where
 $\pr_{\L^{a+1,b}_\circ E}$ denotes the isotypical projection to the
 Cartan summand $\L^{a+1,b}_\circ E\,\subset\,\L^{a+1}E\otimes\L^bE$
 in the tensor product of $\L^{a+1}E$ and $\L^bE$. Evidently the argument
 hinges on the statement that the operator $e\wedge e\wedge e^\flat\ins
 \otimes e\wedge e^\flat\ins\,=\,0$ is trivial; replacing it by the
 analogous statements $e^\flat\ins e\wedge e^\flat\ins\otimes e\wedge
 e^\flat\ins\,=\,0$ etc.~the other seven cases are extremely similar
 and pose no further difficulties.
 \qed

\begin{thebibliography}{22}
%
{\footnotesize
%
\bibitem{besse}
 {\sc A.~Besse},
 {\sl Einstein manifolds},
 Ergebnisse der Mathematik und ihrer Grenzgebiete {\bf 10},
 Springer--Verlag, Berlin (1987).
%
\bibitem{cgh}
 {\sc D.~Calderbank, P.~Gauduchon \& M.~Herzlich},
 {\sl Refined Kato inequalities and conformal weights in Riemannian geometry},
 J.~Funct.~Anal. {\bf 173}, no.~1 (2000), 214---255.
%
\bibitem{galaev}
 {\sc A.~Galaev},
 {\sl Decomposition of the covariant derivative of the curvature
 tensor of a pseudo--K\"ahlerian manifold},
 Ann.~Global Anal.~Geom. {\bf 51}, no.~3 (2017), 245---265.
%
\bibitem{gauduchon}
 {\sc P.~Gauduchon},
 {\sl Structures de Weyl et theoremes d'annulation sur une variete
  conforme autoduale},
  Ann.~Scuola Norm.~Sup.~Pisa Cl.~Sci. {\bf 18}, no.~4 (1991), 563---629.
%
\bibitem{gray}
 {\sc A.~Gray},
 {\sl Compact K\"ahler manifolds with nonnegative sectional curvature},
  Invent.~Math. {\bf 41} (1977), 33---43.
%
\bibitem{aku}
 {\sc K.~Heil, A.~Moroianu \& U.~Semmelmann},
 {\sl Killing and conformal Killing tensors},
 J.~Geom.~Phys. {\bf 106} (2016), 383---400.
%
\bibitem{homma}
 {\sc Y.~Homma},
 {\sl Casimir elements and Bochner identities on Riemannian manifolds},
 Prog.~Math.~Phys. {\bf 34}, Birkh\"auser Boston, Boston (2004).
%
\bibitem{homma1}
 {\sc Y.~Homma},
 {\sl Twisted Dirac operators and generalized gradients},
 Ann.~Global Anal.~Geom.~{\bf 50}, no.~2 (2016), 101---127. 
%
\bibitem{ivanov}
 {\sc S.~Ivanov},
 {\sl Geometry of quaternionic K\"ahler connections with torsion},
 J.~Geom.~Phys.~{\bf 41}, no.~3 (2002), 235--–257. 
%
\bibitem{joyce}
 {\sc D.~Joyce},
 {\sl Compact manifolds with special holonomy},
 Oxford Mathematical Monographs. Oxford University Press, Oxford (2000).
 %
\bibitem{lichnerowicz} 
 {\sc A.~Lichnerowicz},
 {\sl Propagateurs et commutateurs en relativite generale},
 Inst.~Hautes \'Etudes Sci. Publ.~Math. {\bf 10} (1961).
%
\bibitem{au10}
 {\sc A.~Moroianu \& U.~Semmelmann},
 {\sl The Hermitian Laplace operator on nearly K\"ahler manifolds}, 
 Comm.~Math.~Phys. {\bf 294}, no.~1 (2010), 251---272.
%
\bibitem{salamon}
 {\sc S.~Salamon},
 {\sl Riemannian geometry and holonomy groups},
 Pitman Research Notes in Mathematics Series 201,
 Longman Scientific \& Technical, Harlow and John Wiley \& Sons,
 New York (1989).
%
\bibitem{crelle}
 {\sc U.~Semmelmann \& G.~Weingart},
 {\sl Vanishing theorems for quaternionic K\"ahler manifolds},
 J.~Reine Angew.~Math. {\bf 544} (2002), 111---132. 
%
\bibitem{compositio}
 {\sc U. Semmelmann \& G.~Weingart},
 {\sl The Weitzenb\"ock machine},
 Compos.~Math. {\bf 146}, no.~2 (2010), 507---540.
%
\bibitem{verbitsky}
 {\sc M.~Verbitsky},
 {\sl Manifolds with parallel differential forms and K\"ahler
 identities for $\G_2$--manifolds},
 J.~Geom.~Phys. {\bf 61}, no.~6 (2011), 1001---1016.
%
\bibitem{gregor}
 {\sc G.~Weingart},
 {\sl Differential forms on quaternionic K\"ahler manifolds}, 
 Handbook of Pseudo--Riemannian Geometry and Supersymmetry,
 IRMA Lect.~Math.~Theor.~Phys. {\bf 16}, Eur.~Math.~Soc.,
 Zürich (2010), 15---37.
%
\bibitem{secs}
 {\sc G.~Weingart},
 {\sl Moments of sectional curvatures},
 https://arxiv.org/abs/1707.06369 (2017).
%
}
%
\end{thebibliography}
\end{document}